\documentclass[a4paper,11pt,leqno]{article}
\usepackage{amsfonts}
\usepackage{amssymb}
\usepackage{amsmath}
\usepackage{graphicx}

\usepackage{amsthm}
\usepackage{cite}

\def\figurename{Figure}
\makeatletter
\renewcommand{\fnum@figure}[1]{\figurename~\thefigure}

\newtheorem{theorem}{Theorem}[section]

\numberwithin{equation}{section}

\begin{document}

\title{\textbf{Two Non-Congruent Regular Polygons Having Vertices at the Same Distances from the Point}}
\author{\textbf{Mamuka Meskhishvili}}

\date{}
\maketitle

\begin{abstract}
For the given regular plane polygon and an arbitrary point in the plane of the polygon, the distances from the point to the vertices of the polygon are defined. We proved that there is one more non-congruent re\-gu\-lar polygon having the vertices at the same distances from the point. The sizes of both regular polygons are uniquely determined by these distances. In general case, geometrical construction of the second regular polygon is given. It is proved that there are two points in the plane, which separately have the same set of the distances to the vertices of two non-congruent regular polygons with a shared vertex.


Pompeiu theorem, Equilateral triangle
\end{abstract}

\bigskip
\bigskip

\section{Introduction}
\label{sec:1}

\bigskip

The concept of the cyclic averages are introduced in \cite{5}, \cite{6}. For a regular polygon with $n$ vertices $P_n$, there are an $n-1$ number of the cyclic averages:
$$  S_n^{(2)},S_n^{(4)},\dots,S_n^{(2n-2)}.      $$

For an arbitrary point $M$ in the plane of the regular polygon $P_n$, we use the notation $M(d_1,d_2,\dots,d_n,L)$ where $d_i$ are the distances from this point to the vertices $A_i$ of the regular polygon $P_n$ and $L$ is the distance between the point and the center $O$ of the polygon. If the radius of the circumcircle $\Omega$ of the regular polygon $P_n$ is $R$, we denote such polygon by $P_n(R)$.

The cyclic averages are defined as sums of the like even powers of distances $d_i$ to the vertices of $P_n(R)$:
\begin{gather*}
    S_n^{(2)}=\frac{1}{n}\sum_{i=1}^n d_i^2,\;S_n^{(4)}=\frac{1}{n}\sum_{i=1}^n d_i^4,\;\ldots,\;S_n^{(2m)}=\frac{1}{n}\sum_{i=1}^n d_i^{2m}, \\
    \text{where}\;\; m=1,2,\dots,n-1.
\end{gather*}
The cyclic averages can be expressed only in terms of $R$ and $L$:
\begin{equation}\label{eq:1}
    S_n^{(2m)}=(R^2+L^2)^m+\sum_{k=1}^{\lfloor\frac{m}{2}\rfloor} \binom{m}{2k}\binom{2k}{k} R^{2k}L^{2k}(R^2+L^2)^{m-2k}.
\end{equation}

If we are given no more than the distances $d_1,d_2,\dots,d_n$ from a point to the vertices of an $n$-gon, there obviously is an infinity of $n$-gons determined by the $n$ distances. If, however, the polygon is required to be regular, can the $n$ distances uniquely determine the sizes of the polygon? In the present article we investigate this problem.

\bigskip
\bigskip

\section{General case. Existence}
\label{sec:2}

\bigskip

Let us fix the $n$ distances
$$  d_1,d_2,\dots,d_n        $$
and consider the $L$ and $R$ as unknowns.

The number of the cyclic averages is characteristic of the regular polygon but each regular polygon has at least two cyclic averages -- $S_n^{(2)}$ and $S_n^{(4)}$. From \eqref{eq:1} they equal:
\begin{align}
    S_n^{(2)} & =R^2+L^2, \label{eq:2} \\
    S_n^{(4)} & =(R^2+L^2)^2+2R^2L^2. \label{eq:3}
\end{align}

Substituting \eqref{eq:2} into \eqref{eq:3}:
\begin{equation}\label{eq:4}
    2R^2L^2=S_n^{(4)}-(S_n^{(2)})^2.
\end{equation}

The relations \eqref{eq:1}, \eqref{eq:2} and \eqref{eq:4} give us the conditions, which must be satisfied by the $d_1,d_2,\dots,d_n$ if they serve as the distances from the point to the vertices of the regular polygon
$$  S_n^{(2m)}=(S_n^{(2)})^m+\sum_{k=1}^{\lfloor\frac{m}{2}\rfloor} \frac{1}{2^k}\,\binom{m}{2k}\binom{2k}{k} \big(S_n^{(4)}-(S_n^{(2)})^2\big)^k(S_n^{(2)})^{m-2k},      $$
where $m=3,\dots,n-1$. But we initially assumed they are such distances, so we consider \eqref{eq:2} and \eqref{eq:3} only.

From \eqref{eq:2} and \eqref{eq:4} $R^2$ and $L^2$ are the solutions of the equation
$$  X^2-S_n^{(2)}X+\frac{1}{2}\,\big(S_n^{(4)}-(S_n^{(2)})^2\big)=0,        $$
so we get two pairs of the solutions:
\allowdisplaybreaks
\begin{align*}
    \text{I.} \qquad & R_1^2=\frac{1}{2}\,\Big(S_n^{(2)}+\sqrt{3(S_n^{(2)})^2-2S_n^{(4)}}\Big), \\
    & L_1^2=\frac{1}{2}\,\Big(S_n^{(2)}-\sqrt{3(S_n^{(2)})^2-2S_n^{(4)}}\Big); \\[0.2cm]
    \text{II.} \qquad & R_2^2=\frac{1}{2}\,\Big(S_n^{(2)}-\sqrt{3(S_n^{(2)})^2-2S_n^{(4)}}\Big), \\
    & L_2^2=\frac{1}{2}\,\Big(S_n^{(2)}+\sqrt{3(S_n^{(2)})^2-2S_n^{(4)}}\Big).
\end{align*}
If one of them exists, automatically exists another one. Algebraically, it means the following inequalities must be held:
\begin{align}
    3(S_n^{(2)})^2-2S_n^{(4)} & \geq 0, \tag{$*$} \label{eq:*} \\
    S_n^{(2)}-\sqrt{3(S_n^{(2)})^2-2S_n^{(4)}} & \geq 0.    \tag{$**$} \label{eq:**}
\end{align}
Indeed,
$$  3(S_n^{(2)})^2-2S_n^{(4)}=3(R^2+L^2)^2-2\big((R^2+L^2)^2+2R^2L^2\big)=(R^2-L^2)^2,      $$
and from \eqref{eq:4} follows
$$  S_n^{(4)}\geq (S_n^{(2)})^2,        $$
which proves \eqref{eq:**}.

Denote by $\Omega_1$ and $\Omega_2$ the circumcircles of the regular polygons $P_n(R_1)$ and $P_n(R_2)$, by $O_1$ and $O_2$ their centers, respectively. Therefore
$$  L_1=MO_1 \;\;\text{and}\;\; L_2=MO_2.       $$
From the solutions I and II follows:
$$  R_1>R_2 \;\;\text{and}\;\; L_1<L_2,     $$
so the first solution corresponds to the larger regular polygon, while the second solution corresponds to the smaller one. The distances from the $M$ point to the centers is longer for the smaller polygon.

If
$$  3(S_n^{(2)})^2=2S_n^{(4)},       $$
the point $M$ lies on the circumcircle,
$$  R_1=L_1=R_2=L_2.     $$
This is degenerate case -- both regular polygons are congruent. We obtain:

\begin{theorem}\label{th:1}
If the point of the distances $d_1,d_2,\dots,d_n$ to the vertices of the regular polygon $P_n(R_1)$, does not lie on the circumcircle $\Omega_1$ of $P_n(R_1)$, there is one more non-congruent regular polygon $P_n(R_2)$ having the vertices at the same $d_1,d_2,\dots,d_n$ distances from the point. If the point lies on the circumcircle $\Omega_1$ of $P_n(R_1)$ there is no more non-congruent regular polygon having the vertices at the same distances from the point.
\end{theorem}

For the first solution (larger polygon) $L_1<R_1$ i.e. the point $M$ lies inside the circumcircle $\Omega_1$, while for the second solution (smaller polygon) $L_2>R_2$ i.e point $M$ lies outside the circumcircle $\Omega_2$.

Let us summarize the obtained results.

\begin{theorem}\label{th:2}
If the arbitrary point $M(d_1,d_2,\dots,d_n,L_1)$ lies inside the circumcircle $\Omega_1$ of the regular polygon $P_n(R_1)$, there is one more regular polygon $P_n(R_2)$ having the vertices at the same distances from the point $M(d_1,d_2,\dots,d_n,L_2)$ and holds:
\allowdisplaybreaks[0]
\begin{multline*}
    \sqrt{\frac{1}{2}\,\Big(S_n^{(2)}+\sqrt{3(S_n^{(2)})^2-2S_n^{(4)}}\Big)}=R_1=L_2 \\
    >R_2=L_1=\sqrt{\frac{1}{2}\,\Big(S_n^{(2)}-\sqrt{3(S_n^{(2)})^2-2S_n^{(4)}}\Big)}\,,
\end{multline*}
i.e. the point $M$ lies outside the circumcircle $\Omega_2$ of the regular polygon $P_n(R_2)$.
\end{theorem}

\begin{theorem}\label{th:3}
If the arbitrary point $M(d_1,d_2,\dots,d_n,L_1)$ lies outside the circumcircle $\Omega_1$ of the regular polygon $P_n(R_1)$, there is one more regular polygon $P_n(R_2)$ having the vertices at the same distances from the point $M(d_1,d_2,\dots,d_n,L_2)$ and holds:
\begin{multline*}
    \sqrt{\frac{1}{2}\,\Big(S_n^{(2)}-\sqrt{3(S_n^{(2)})^2-2S_n^{(4)}}\Big)}=R_1=L_2 \\
    <R_2=L_1=\sqrt{\frac{1}{2}\,\Big(S_n^{(2)}+\sqrt{3(S_n^{(2)})^2-2S_n^{(4)}}\Big)}\,,
\end{multline*}
i.e. the point $M$ lies inside the circumcircle $\Omega_2$ of the regular po\-ly\-gon $P_n(R_2)$.
\end{theorem}

\newpage
\bigskip
\bigskip

\section{Equilateral triangle}
\label{sec:3}

\bigskip
\noindent \textit{\noindent Algebraic Backround}
\medskip

The well-known the Pompeiu theorem states \cite{8}:

let given an equilateral triangle and any point in its plane. Then the distances from the point to the vertices $d_1$, $d_2$, $d_3$ are lengths of the sides of a triangle.

we call a triangle with sides $d_1$, $d_2$, $d_3$ a Pompeiu triangle \cite{9}. The Pompeiu triangle is degenerate if the point lies on the circumcircle of the equilateral triangle, because by Van Schooten's theorem the largest distance equals to the sum of the others.

According to Theorem \ref{th:1}, for the given Pompeiu triangle there are two equilateral triangles -- the larger and the smaller. In \cite{1}, \cite{3}, \cite{4}, \cite{9} are investigated case of the larger equilateral triangle i.e. the point lies inside of the circumcircle, both equilateral triangles are considered by H.~Eves \cite{7}.

For the equilateral triangle $P_3(R)$ and the point $M(d_1,d_2,d_3,L)$ the $S_3^{(2)}$ and $S_3^{(4)}$ cyclic averages equal:
$$  S_3^{(2)}=\frac{1}{3}\,(d_1^2+d_2^2+d_3^2), \quad S_3^{(4)}=\frac{1}{3}\,(d_1^4+d_2^4+d_3^4).   $$
The inequality \eqref{eq:*} gives
\begin{align*}
    3(S_n^{(2)})^2-2S_n^{(4)} & =\frac{1}{3}\,\big((d_1^2+d_2^2+d_3^2)^2-2(d_1^4+d_2^4+d_3^4)\big) \\
    & =\frac{16}{3}\,\Delta^2(d_1,d_2,d_3),
\end{align*}
the symbol -- $\Delta_{(d_1,d_2,d_3)}$ denotes the area of the triangle whose sides have lengths $d_1$, $d_2$, $d_3$. Then the inequality \eqref{eq:**} turns into will-known   \linebreak   Weit\-zen\-b\"{o}ck's inequality \cite{2}
$$  d_1^2+d_2^2+d_3^2\geq 4\sqrt{3}\,\Delta_{(d_1,d_2,d_3)}.        $$

Two equilateral triangles are:
\begin{enumerate}
\item[I.] The larger triangle, $M$ lies inside $\Omega_1$;
\begin{align*}
    R_1^2 & =\frac{1}{6}\,\Big(d_1^2+d_2^2+d_3^2+4\sqrt{3}\,\Delta_{(d_1,d_2,d_3)}\Big), \\
    L_1^2 & =\frac{1}{6}\,\Big(d_1^2+d_2^2+d_3^2-4\sqrt{3}\,\Delta_{(d_1,d_2,d_3)}\Big).
\end{align*}

\item[II.] The smaller triangle, $M$ lies outside $\Omega_2$;
\begin{align*}
    R_2^2 & =\frac{1}{6}\,\Big(d_1^2+d_2^2+d_3^2-4\sqrt{3}\,\Delta_{(d_1,d_2,d_3)}\Big), \\
    L_2^2 & =\frac{1}{6}\,\Big(d_1^2+d_2^2+d_3^2+4\sqrt{3}\,\Delta_{(d_1,d_2,d_3)}\Big).
\end{align*}
\end{enumerate}

\medskip
\noindent \textit{Geometrical Construction}
\medskip

We perform the constructions in two ways:
\begin{enumerate}
\smallskip
\item[A.] Given the Pompeiu triangle and construct both equilateral triangles.

\smallskip
\item[B.] Given one of the equilateral triangle and the point and construct the second one.
\end{enumerate}

\begin{figure}[h]
\centerline{\includegraphics[width=8cm]
    {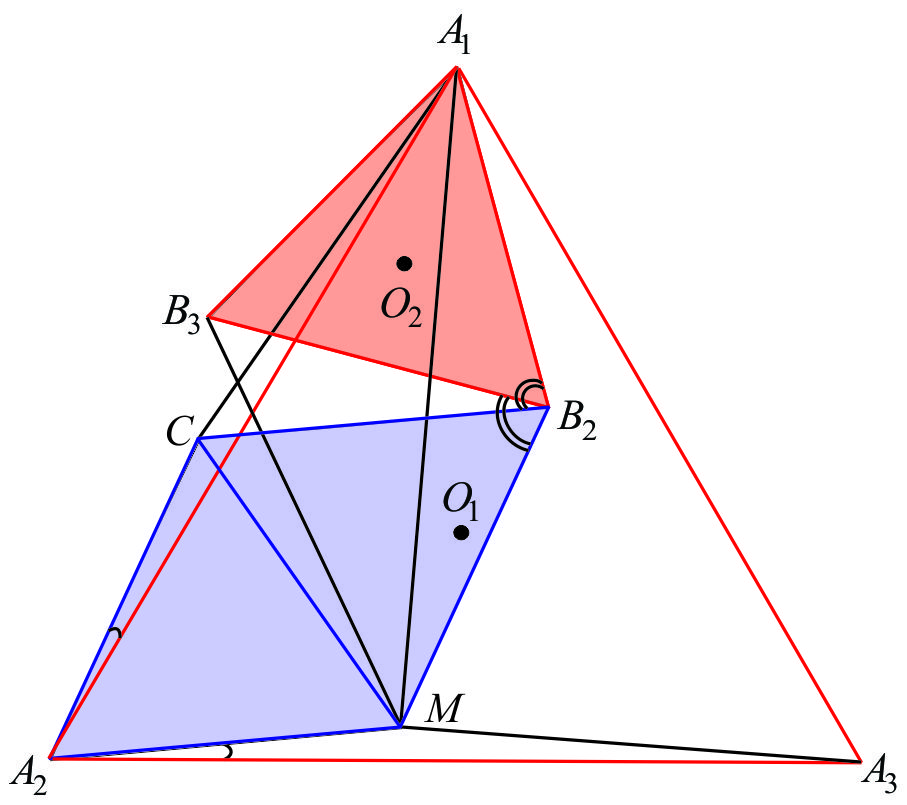} }
\caption{}
\label{fig:1}
\end{figure}

\begin{figure}[h]
\centerline{\includegraphics[width=8.7cm]
    {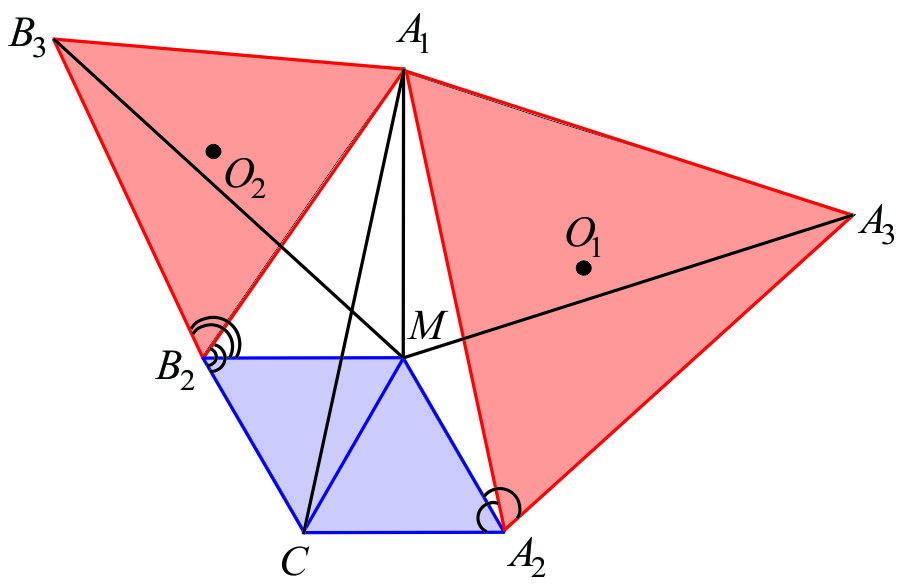} }
\caption{}
\label{fig:2}
\end{figure}

%\bigskip
\medskip
A. If given the Pompeiu triangle we take one vertex as the point $M$, see Fig.~\ref{fig:1} and Fig.~\ref{fig:2}. If $MCA_1$ is the Pompeiu triangle we construct around one side, for example $MC$, two auxiliary equilateral triangles $MCA_2$ and $MCB_2$. Obtained the auxiliary points $A_2$ and $B_2$ connect to the third vertex $A_1$. Two line segments $A_2A_1$ and $B_2A_1$ serve as the sides of the desired two equilateral triangles -- $A_1A_2A_3$ and $A_1B_2B_3$. Indeed,
$$  MA_3=CA_1=MB_3;     $$
because of the congruency of the triangles:
$$  MA_2A_3=CA_2A_1 \;\;\text{and}\;\; CB_2A_1=MB_2B_3.     $$

\begin{figure}[h]
\centerline{\includegraphics[width=9cm]
    {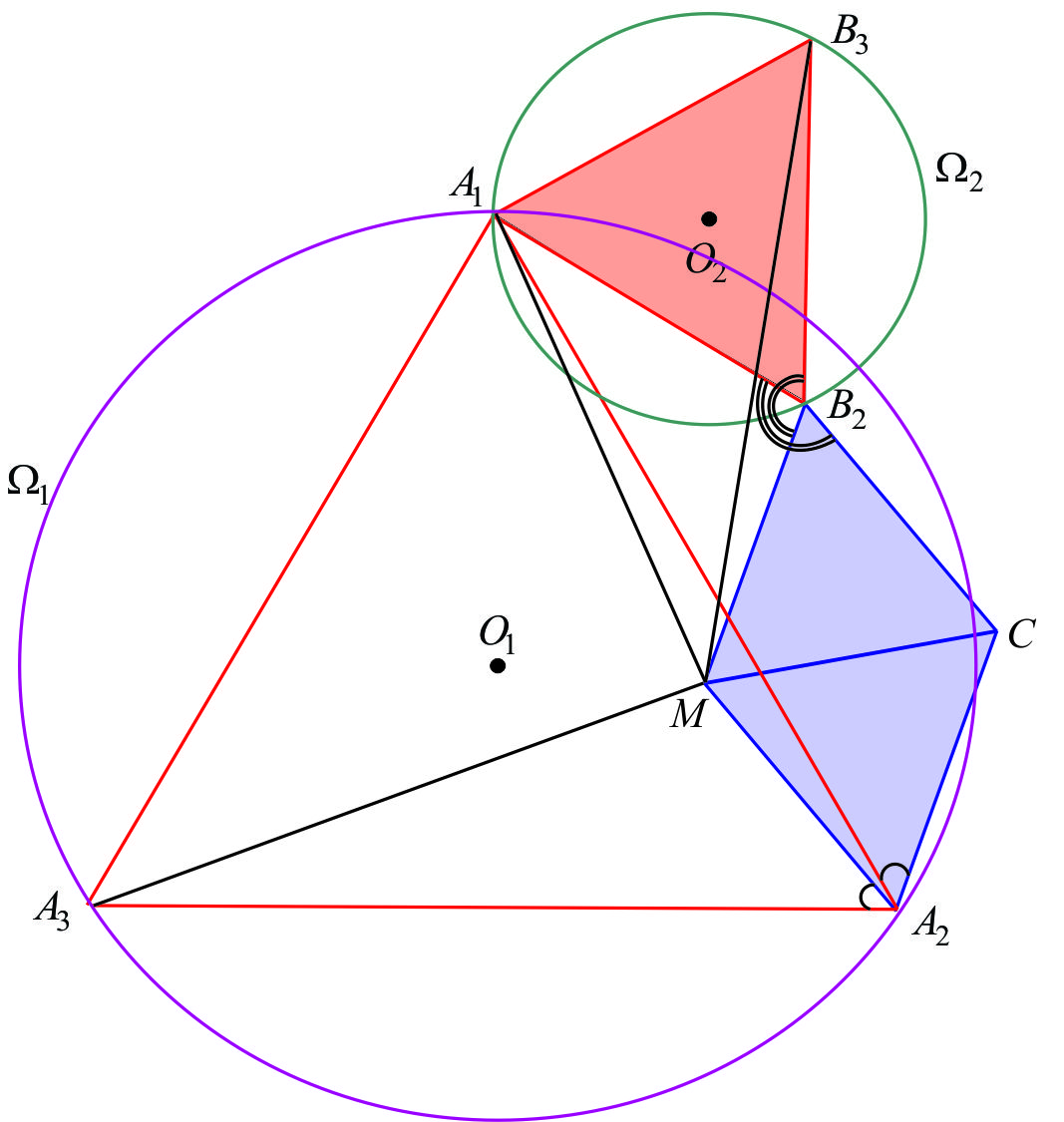} }
\caption{}
\label{fig:3}
\end{figure}

B. If the point $M$ inside of the circumcircle $\Omega_1$, it means the triangle $A_1A_2A_3$ is the larger, see Fig.~\ref{fig:3}. We construct around one distance, for example $MA_2$, one auxiliary equilateral triangle $MA_2C$ obtained the auxiliary point $C$. Again around line segment $MC$ construct the second auxiliary equilateral triangle -- $MCB_2$. Connect the obtained point $B_2$ to the vertex $A_1$. The line segment $B_2A_1$ serve as the side of the desired smaller equila\-te\-ral triangle $A_1B_2B_3$.

\begin{figure}[h]
\centerline{\includegraphics[width=7cm]
    {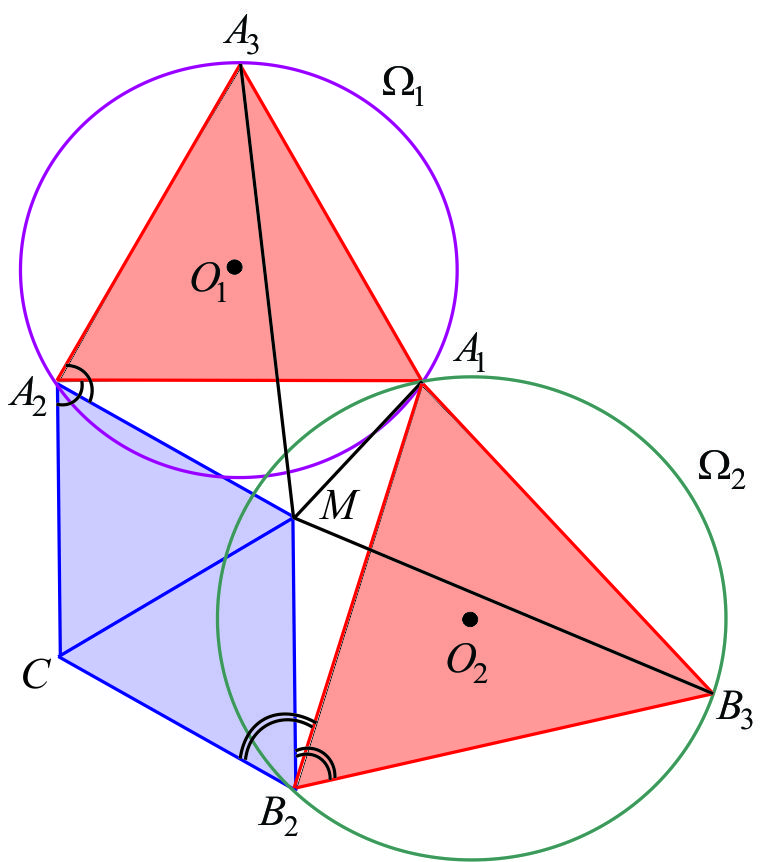} }
\caption{}
\label{fig:4}
\end{figure}

If the point $M$ outside of the circumcircle $\Omega_1$, it means the triangle $A_1A_2A_3$ is the smaller, see Fig.~\ref{fig:4}. Repeat above-mentioned steps, we obtain the larger equilateral triangle -- $A_1B_2B_3$.

%\bigskip
%\medskip

%\bigskip
%\medskip

\bigskip
\bigskip

\section{General case. Construction}
\label{sec:4}

\bigskip

For the regular polygon $P_n$, when $n>3$ constructions by using the auxiliary triangles are impossible. In general case Theorem \ref{th:2} and Theorem \ref{th:3} give us the method of construction for the second regular polygon from the given one and the point. From these theorems:
\begin{equation}\label{eq:***}
    R_2=L_1 \;\;\text{and}\;\; L_2=R_1.  \tag{$*\!*\!*$}
\end{equation}
Conditions \eqref{eq:***} are necessary but not sufficient for the construction.

Let us consider the distances from the given point $M$ to the vertices of the second regular polygon as unknowns:
$$  x_1,x_2,\dots,x_n.      $$

From \eqref{eq:1} and \eqref{eq:***} they satisfy
\begin{equation}\label{eq:4.1}
    \sum_{i=1}^n x_i^{2m}=\sum_{i=1}^n d_i^{2m}, \;\;\text{where}\;\; m=1,2,\dots,n-1.
\end{equation}

We have the $n-1$ equations, so to determine the unknowns uniquely, let us consider one unknown as one of the distances -- $d_i$. Without loss of generality take
\begin{equation}\label{eq:4.2}
    x_1=d_1.
\end{equation}
Then, from the elementary properties of the symmetric functions, from \eqref{eq:4.1} it follows that
$$  d_2^2,\dots,d_n^2 \;\;\text{and}\;\; x_2^2,\dots,x_n^2      $$
are roots of the same equation of degree $n-1$. Consequently $x_2,\dots,x_n$ are a permutation of the $d_2,\dots,d_n$; therefore both of them are the same set of the distances:
$$  \big\{x_2,\dots,x_n\big\}=\big\{d_2,\dots,d_n\big\}.      $$
So, the conditions \eqref{eq:***} and \eqref{eq:4.2} are the sufficient conditions to identify the second regular polygon having the vertices at the same distances from the point.

\begin{figure}[ht]
\centerline{\includegraphics[width=8cm]     %6.5cm]
    {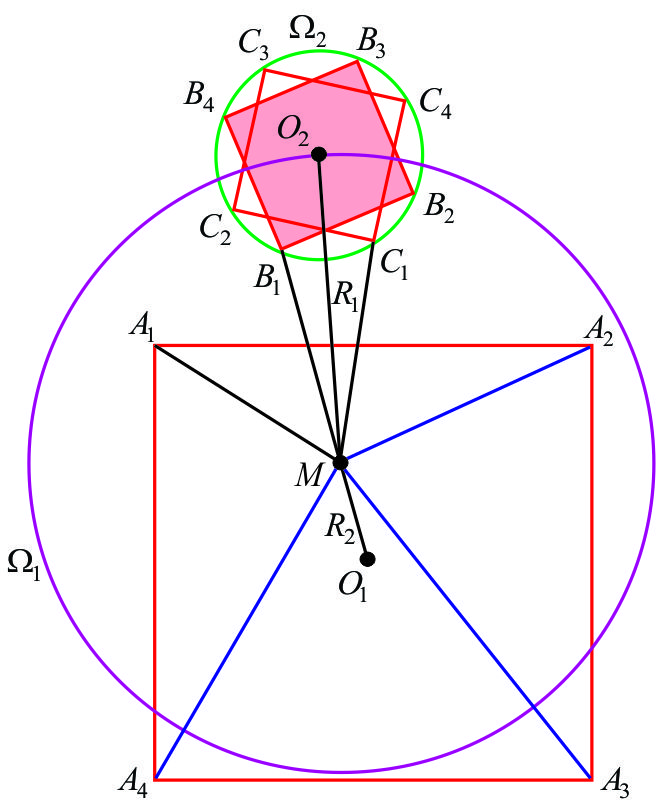} }
\caption{}
\label{fig:5}
\end{figure}

\begin{figure}[h]
\centerline{\includegraphics[width=9cm]     %7.5cm]
    {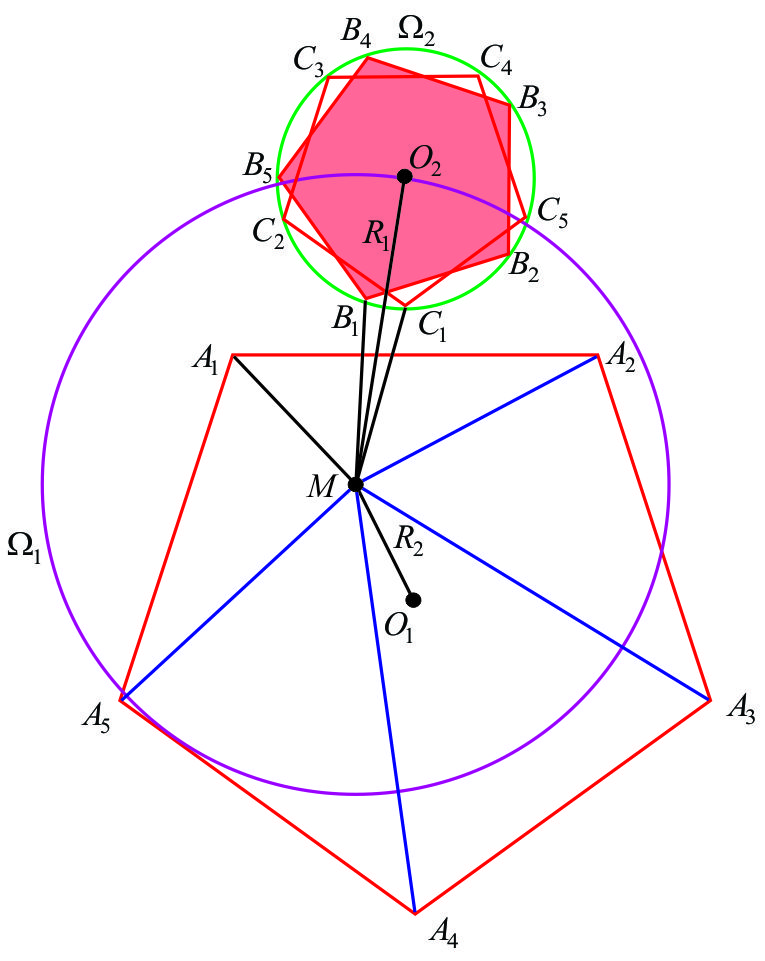} }
\caption{}
\label{fig:6}
\end{figure}

Constructions for a square and a regular pentagon are given in Fig.~\ref{fig:5} and Fig.~\ref{fig:6}, respectively. We describe the method of construction, which is true for any regular polygon.

%\bigskip

If given regular polygon $A_1A_2\cdots A_n$ and the point $M$ in its plane, draw the circle $\Omega_1$ with center $M$ and the radius $R_1$. Choose the point $O_2$ on the $\Omega_1$, which is the center of the second desired regular polygon -- $B_1B_2\cdots B_n$. Construct the circle $\Omega_2$, whose center is $O_2$ and the radius equals $R_2=O_1M$. From the point $M$ as the center draw the auxiliary circle with radius $MA_1$. The intersection points of the auxiliary and $\Omega_2$ circles -- $B_1$ and $C_1$ are the vertices of the desired polygon separately. In fact, we construct two regular polygons -- $B_1B_2\cdots B_n$ and $C_1C_2\cdots C_n$ which have the same set of the distances from the point $M$ to the vertices of original $A_1A_2\cdots A_n$ polygon:
\begin{gather*}
    \big\{MA_i\big\}=\big\{MB_i\big\}=\big\{MC_i\big\}, \\
    \text{where}\;\; i=1,2,\dots,n.
\end{gather*}
For the square case, see Fig.~\ref{fig:5}, by corresponding enumeration of the vertices:
\begin{gather*}
    MA_i=MB_i=MC_i, \;\;\text{where}\;\; i=1,\dots,4.
\end{gather*}
For the regular pentagon case, see Fig.~\ref{fig:6}:
\begin{gather*}
    MA_i=MB_i=MC_i, \;\;\text{where}\;\; i=1,\dots,5.
\end{gather*}

\bigskip
\bigskip

\section{Two regular polygons and two points theorem}
\label{sec:5}

\bigskip

Until now we consider the case when the regular polygon and the point were given. Let us consider the case when two non-congruent regular polygons (with the same number of vertices) are given initially. Is there a point in the plane of two non-congruent regular polygons, from where the distances to the vertices of these polygons are the same? From the construction method it is clear -- such point exists, if
$$  |R_1-R_2|\leq O_1O_2\leq R_1+R_2,       $$
and one of the distances $MA_i$ and $MB_i$ must be equal to each other \eqref{eq:4.2}.
If two polygons have one shared vertex this equality automatically holds. From the construction method the point $M$ must be at distance $R_1$ from the $O_2$, and at distance $R_2$ from the $O_1$, i.e. intersection of two circles -- $\Omega_1(O_2,R_1)$ and $\Omega_2(O_1,R_2)$. But in general case, there are two points of such properties. So we obtain -- two regular polygons and two points theorem.

\bigskip

\begin{figure}[h]
\centerline{\includegraphics[width=10cm]     %8.5cm]
    {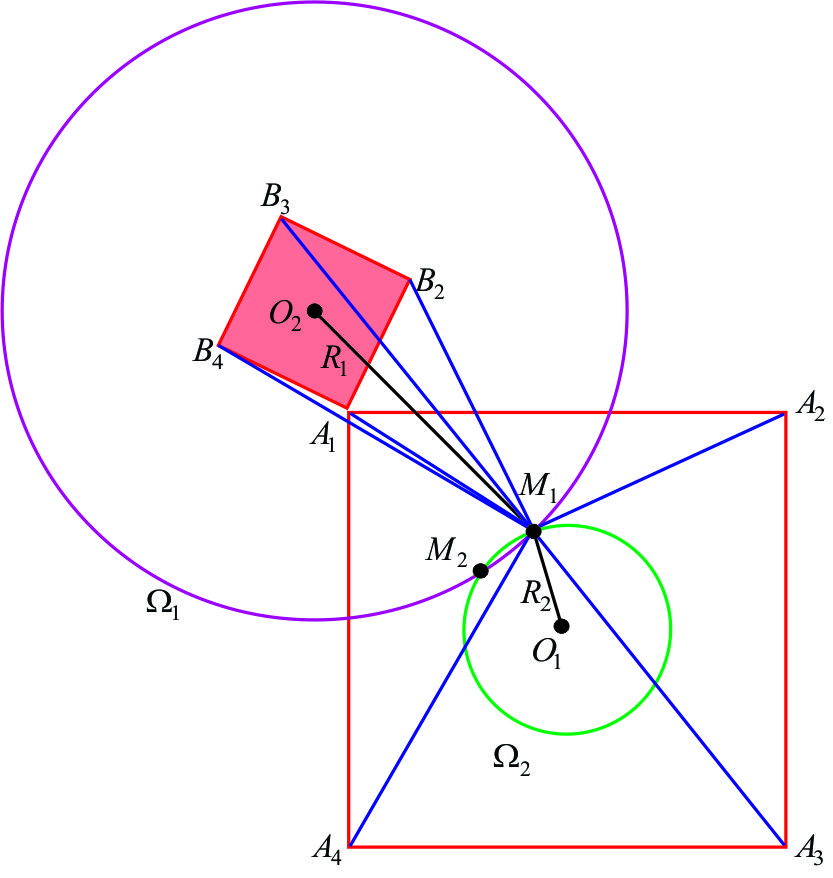} }
\caption{}
\label{fig:7}
\end{figure}

\bigskip

\begin{theorem}\label{th:4.1}
If two non-congruent regular polygons $A_1A_2\cdots A_n$ and    \linebreak    $B_1B_2\cdots B_n$ have one shared vertex, there are two points $M_1$ and $M_2$ in the plane, separately having the same set of the distances to the vertices of the polygons:
\begin{gather*}
    \big\{M_1A_i\big\}=\big\{M_1B_i\big\}, \quad \big\{M_2A_i\big\}=\big\{M_2B_i\big\}, \quad \text{where}\;\; i=1,\dots,n.
\end{gather*}
The $M_1$ and $M_2$ are intersection points of two circles -- $\Omega_1(O_2,R_1)$ and $\Omega_2(O_1;R_2)$, where $O_1$ and $O_2$ are the centers of the circumcircles of   \linebreak    $A_1A_2\cdots A_n$ and $B_1B_2\cdots B_n$, $R_1$ and $R_2$ their radii, respectively.
\end{theorem}

The construction of $M_1$ and $M_2$ is given for the squares in Fig.~\ref{fig:7}, the shared vertex is $A_1$.

The same distances are:
$$  M_1A_2=M_1B_2, \quad M_1A_3=M_1B_3, \quad M_1A_4=M_1B_4     $$
and
$$  M_2A_2=M_2B_4, \quad M_2A_3=M_2B_3, \quad M_2A_4=M_2B_2.    $$
If the shared vertex and the centers of the circumcircles are collinear, there is only one point having the same set of the distances to the vertices of two non-congruent regular polygons.

\bigskip

\bigskip

\bigskip

DEPARTMENT OF MATHEMATICS

GEORGIAN-AMERICAN HIGH SCHOOL

18 CHKONDIDELI STR., TBILISI 0180, GEORGIA

\textit{E-mail address:} \texttt{mathmamuka@gmail.com}

\end{document}